\documentclass[12pt,leqno,fleqn]{amsart}  
\usepackage{amsmath,amstext,amsthm,amssymb,amsxtra}
\usepackage{txfonts,pxfonts} 
\usepackage[colorlinks,citecolor=red,pagebackref,hypertexnames=false]{hyperref}
\usepackage[backrefs,author-year]{amsrefs}
\allowdisplaybreaks
\swapnumbers
 
\theoremstyle{plain} 
\newtheorem{lemma}[equation]{Lemma} 
\newtheorem{proposition}[equation]{Proposition} 
\newtheorem{theorem}[equation]{Theorem} 

\newtheorem{cww}[equation]{Chang-Wilson-Wolff Inequality } 
\newtheorem*{halasz}{Hal{\'a}sz' Theorem}
\newtheorem*{roth}{K.~Roth's Theorem}
\newtheorem*{ell1}{Conjecture on the $ L ^{1}$ norm of  $ D_N$} 
\newtheorem*{Lp}{Conjecture on the $ L ^{p}$ norm of  $ D_N$} 
\newtheorem*{Hardy}{An Equivalent Definition of Real Valued Hardy Space}
\newtheorem*{rademacher}{Khintchine Inequalities for Rademacher Random Variables}

\theoremstyle{definition}
\newtheorem{definition}[equation]{Definition} 

\theoremstyle{remark}

\newtheorem*{acknowledgment}{Acknowledgment}

\numberwithin{equation}{section}

\def\norm#1.#2.{\lVert#1\rVert_{#2}}
\def\Norm#1.#2.{\bigl\lVert#1\bigr\rVert_{#2}}
\def\NOrm#1.#2.{\Bigl\lVert#1\Bigr\rVert_{#2}}
\def\NORm#1.#2.{\biggl\lVert#1\biggr\rVert_{#2}}
\def\NORM#1.#2.{\Biggl\lVert#1\Biggr\rVert_{#2}}

\def\ip#1,#2,{\langle #1,#2\rangle}
\def\Ip#1,#2,{\bigl\langle#1,#2\bigr\rangle}
\def\IP#1,#2,{\Bigl\langle#1,#2\Bigr\rangle}

\def\abs#1{\lvert#1\rvert}

\def\ABs#1{\biggl\lvert#1\biggr\rvert}

\def\XXint#1#2#3{{\setbox0=\hbox{$#1{#2#3}{\int}$}
     \vcenter{\hbox{$#2#3$}}\kern-.5\wd0}}

\begin{document}
\title[The Discrepancy Function Close to $ L ^{1}$]{On the Discrepancy Function in Arbitrary Dimension,\\ Close to $ L ^{1}$}
\author{Michael Lacey}
\email{lacey@math.gatech.edu}

\begin{abstract}
Let  $\mathcal A_N$ to be $N$ points in the unit cube in dimension $ d$, and consider
the Discrepancy function 
\begin{equation*} 
D_N( \vec x) \coloneqq \sharp \big(\mathcal A_N \cap [\vec 0,\vec x)\big)
-N \abs{ [\vec 0,\vec x)}
\end{equation*}
Here, $ \vec x= (x_1 ,\dotsc, x_d)$, 
$[ 0,\vec x)=\prod _{t=1} ^{d} [0,x_t)$, and $ \lvert[\vec 0,\vec x)  \rvert $ 
denotes the Lebesgue measure of the rectangle.  
We show that necessarily 
\begin{equation*}
\norm D_N. L ^{1} (\log L) ^{ (d-2)/2}. \gtrsim (\log N) ^{ (d-1)/2}\,. 
\end{equation*}
In dimension $ d=2$, the `$ \log L$' term has power zero, which corresponds to 
a Theorem due to   \cite{MR637361}. The power on $ \log L$ in dimension  $ d\ge 3$ appears 
to be new, and supports a well-known conjecture on the $L ^{1}$ norm of $ D_N$.  Comments 
on the Discrepancy function in Hardy space also support the conjecture. 
\end{abstract}

\maketitle

\section{Main Theorem} 

Our subject is irregularities of distribution of points with respect to 
rectangles in the unit cube.  
It is a familiar theme of the subject is to show that no matter how $ N$ 
points are selected, they must be  far from uniform.  
We give a new proof of a well-known theorem in the subject  \cite{MR637361}, 
concerning the $ L ^{1}$ norm of the Discrepancy function, 
and show that this result admits an extension to arbitrary dimension. 
We also make some remarks on the Discrepancy function and Hardy space. 

Let $ \mathcal A_N \subset [0,1] ^{d}$ be a set of cardinality $ N$.  Define the 
\emph{Discrepancy Function} associated to $ \mathcal A_N$ as 
a function on the unit square as follows. 
\begin{equation*}
D_N (\vec x) \coloneqq \sharp \bigl(\mathcal A_N \cap [0, \vec x)\bigr) - N \abs{ [0, \vec x)}\,,
\end{equation*}
where $ [0,\vec x)$ is the rectangle in the unit cube with one vertex at the origin 
and the other at $ \vec x =(x_1 ,\dotsc, x_d)$, and $  \abs{ [0, \vec x)}$ denotes  
the Lebesgue measure of the rectangle.  
This is the difference between the number of points 
in the rectangles $ [0,\vec x)$ and the expected number of points in the rectangle.
The relative size of this function, in many senses, 
necessarily must increase with $ N$.    The principal result here is that 
of \cite{MR0066435}.  

\begin{roth} We have the estimate below, valid in all dimensions $ d\ge 2$
\begin{equation*}
\norm D_N. 2. \gtrsim (\log N) ^{(d-1)/2}
\end{equation*}
where the implied constant is only a function of dimension $ d$. 
\end{roth}

The same bound holds for the $ L ^{p}$ norm, for $ 1<p<\infty $, \cite{MR0491574}, 
and is known to be sharp as to the order of magnitude. 
The endpoint cases of $ p=1$ and $ p=\infty $ are much harder.  We concentrate on the 
case of $ p=1$ in this note, and refer the reader to 
\cites{MR1032337,0705.4619,math.CA/0609815,MR637361}  for more information about 
the case of $ p=\infty $.

In the case of $ d=2$, we have definitive information about the $ L ^{1}$ norm, 
namely that the Roth lower bound holds.  See \cite{MR637361}. 

\begin{halasz} In dimension $ d=2$ we have the uniform estimate 
\begin{equation*}
\norm D_N .1. \gtrsim \sqrt {\log N}
\end{equation*}
\end{halasz}

A principal conjecture in the subject is that the Roth bound holds for the $ L ^{1}$ 
norm of the Discrepancy function in all dimensions. 

\begin{ell1} In dimension $ d\ge 3$ we have the estimate 
\begin{equation*}
\norm D_N .1. \gtrsim  (\log N ) ^{ (d-1)/2}\,. 
\end{equation*}
\end{ell1}

The best known result for the $ L ^{1}$ norm directly is the Hal{\'a}sz bound below. 
\begin{equation*}
\norm D_N .1. \gtrsim \sqrt {\log N}\,, \qquad d\ge 3\,. 
\end{equation*}
This is a simple consequence of his argument in \cite{MR637361}. 

Our main result is a partial extension of Hal{\'a}sz' Theorem to arbitrary dimension. 

\begin{theorem}\label{t.beyondHalasz} In dimension $ d\ge 2$ we have 
\begin{equation*}
\norm D_N. L ^{1} (\log L) ^{(d-2)/2}. \gtrsim (\log N) ^{ (d-1)/2}\,. 
\end{equation*}
\end{theorem}

Here, we use an Orlicz norm which is `close' to $ L ^{1}$. 
Its definition, the well-known one, is made precise in the next section. 
We remark that the 
orthogonal function method of  \cite{MR0066435}, especially as modified by the 
observations in \cite{MR0491574}, can be used to prove the estimate 
\begin{equation*}
\norm D_N. L ^{1} (\log L) ^{(d-1)/2}. \gtrsim (\log N) ^{ (d-1)/2}\,. 
\end{equation*}
This does not contain Hal{\'a}sz' Theorem, as the power of the $ \log L$ is too large.
Our proof is an elaboration of 
that of Hal{\'a}sz, using appropriate version of the Chang-Wilson-Wolff inequality 
\cite{MR800004}, and the variant in \cite{MR1439553}.

Hal{\'a}sz' proof is by way duality,\footnote{Extremal choices of the point set will 
result in a Distribution Function supported whose $ L^1$ norm is determined on 
a set which is nearly the whole square. A proof by duality is natural. }
namely one constructs an appropriate 
bounded function $ \Psi $, and  obtains a uniform lower bound
on the inner product $\ip \Psi , D_N, $.
And in particular, the proof uses a Bernoulli product construction. 
Our construction is not Bernoulli product, though once the function is constructed, 
many details are variants of the arguments in \cite{MR637361}.  

The concluding section of the paper includes some remarks about the 
Discrepancy function and Hardy spaces, and proves a result which can be 
thought of supporting evidence for the $ L ^{1}$ conjecture above.

\begin{acknowledgment}
The author benefited from many stimulating conversations with Mihalis Kolountzakis on the 
subject of this paper. 
\end{acknowledgment}

\section{Preliminary Facts} 

We suppress many constants which do not affect the arguments in essential ways. 
$ A \lesssim B$ means that there is an absolute constant $ K>0$ so that $ A \le K B$. 
Thus $ A \lesssim 1$ means that $ A$ is bounded by an absolute constant. 
And if $ A \lesssim B \lesssim A$, we write $ A \simeq B$.

\bigskip

\subsection*{Inequalities}

We collect some standard estimates from the Probability and Harmonic Analysis 
literature.  Let $ \{r_j \mid j\ge 1\}$ be a sequence of Rademacher random variables, 
thus independent identically distributed  random variables 
on a probability space $ (\Omega , \mathbb P )$ such that $ \mathbb P ( r_j=1)= \mathbb P
(r_j=-1)=\tfrac 12$.  We have the moment function inequality 
\begin{equation*}
\mathbb E \operatorname {exp} \Biggl(\lambda \sum _{j} c_j r_j \Biggr)
\le \operatorname {exp} \Biggl( \tfrac 12 \lambda ^2 \sum _{j} c_j ^2  \Biggr)\,, 
\qquad \lambda \in \mathbb R \,. 
\end{equation*}
This holds for all sequences of coefficients $ \{c_j\mid j\in \mathbb N \}$ such that 
the right hand side above is finite.  This implies the distributional inequality 
\begin{equation*}
\mathbb P \Biggl( \sum _{j} c_j r_j > t  \Biggr) \lesssim 
\operatorname {exp} \Biggl( -\frac {t ^2 } { 2\sum _{j} c_j ^2 } \Biggr)\,, 
\qquad t>0\,. 
\end{equation*}
An equivalent formulation is in terms of the \emph{Khintchine} inequalities. 

\begin{rademacher} We have the inequalities
\begin{equation*}
\NORm \sum _{j} c_j r_j . p.  \lesssim 
C \sqrt p \Biggl[  \sum _{j} c_j ^2  \Biggr] ^{1/2}  \,,  
\qquad 0<p<\infty \,. 
\end{equation*}
\end{rademacher}

The best constants in these inequalities are of significant interest. For the range 
$ p\ge 2$, of interest to us, see \cite{MR0438089}. 

\bigskip

There is a powerful extension to these inequalities to the setting of 
Haar series, or more generally, conditionally symmetric martingales.  
We state the results in a convenient form.  

In one dimension, the class of dyadic intervals in the unit interval  are 
$\mathcal D {} \coloneqq {}\{ [j2^{-k},(j+1)2^{-k})
\mid j,k\in \mathbb N\,, 0\le j < 2 ^{k}\} $.  Let $ \mathcal D_n$ denote the 
dyadic intervals of length $ 2 ^{-n}$, and by abuse of notation, also the 
sigma field generated by these intervals.  For an integrable function $ f$ on 
$ [0,1]$, the conditional expectation is 
\begin{equation*}
f_n=\mathbb E (f \mid \mathcal D_n) \coloneqq  \sum _{I\in \mathcal D_n} \mathbf 1_{I}  \cdot 
\lvert  I\rvert ^{-1}  \int _{I} f (y)\;dy\,. 
\end{equation*}
The sequence of functions $ \{ f_n \mid n\ge 0\}$ is a \emph{martingale}.  The 
\emph{martingale difference sequence } is $ d_0=f_0$, and $ d_n= f_n-f _{n-1}$ for 
$ n\ge 1$.  The sequence of functions $ \{d_n\mid n\ge 0\}$ are pairwise orthogonal. 
The \emph{square function} is 
\begin{equation*}
\operatorname S (f) \coloneqq \Biggl[\sum _{n=0} ^{\infty } \lvert  d_n\rvert ^2   \Biggr]
^{1/2} \,. 
\end{equation*}
We have the following extension of the Khintchine inequalities.  

\begin{theorem}\label{t.martKhintchine} The inequalities below hold, for 
some absolute choice of constant $ C>0$. 
\begin{equation}\label{e.martKhintchine}
\norm f. p. \le C \sqrt p \norm \operatorname S (f).p.\,, \qquad 2\le p < \infty \,. 
\end{equation}
In particular, this inequality holds for Hilbert space valued functions $ f$. 
\end{theorem}

For real-valued martingales, this was observed by \cite{MR800004}.  The extension 
to Hilbert space valued martingales is useful for us, and is proved in \cite{MR1439553}. 
Indeed, the best constants in these inequalities are known for $ p\ge 3$ \cite{MR1018577}.

\subsection*{Orlicz Spaces}

For background on Orlicz Spaces, we refer the reader to \cite{MR0500056}. 
Consider a symmetric convex function $ \psi $, 
which is zero at the origin,  and is otherwise non-zero.  Let $ (\Omega , P)$ 
be a probability space, on which our functions are defined, and let  $ \mathbb E $ 
denote expectation over the probability space. 
We can define 
\begin{equation}\label{e.psi}
\norm f. L ^\psi . = \inf \{ K>0\mid  \mathbb E \psi (f \cdot K ^{-1} )\le 1\}\,, 
\end{equation}
where we define the infimum over the empty set to be $ \infty $.  The 
set of functions $ L ^{\psi } = \{f \mid  \norm f. L ^{\Psi } . < \infty \}$ 
is a normed linear space, with norm as above.  It is the Orlicz space associated with 
$ \psi $.  

We are interested in, for instance, $ \psi (x)= \operatorname e ^{x ^2 }-1$, in 
which case we indicate the Orlicz space as $ \operatorname {exp} (L ^{2})$.  
More generally, for $ 0<\alpha <1$, we let $ \psi _{\alpha } (x)$ be a symmetric 
convex function which equals $ \operatorname e ^{\lvert  x\rvert ^{\alpha } }-1$ 
for $ \lvert  x\rvert $ sufficiently large, depending upon $ \alpha $.\footnote{We are only interested in 
measuring the behavior of functions for large values of $ f$, so this requirement 
is sufficient. For $ \alpha >1$, we can 
insist upon this equality for all $ x$.}  And we write $ L ^{\psi _{\alpha }} = 
\operatorname {exp}(L ^{\alpha })$.

The following proposition is well-known, and follows from elementary methods. 

\begin{proposition}\label{p.comparable} We have the following equivalence of norms 
valid for all $ \alpha >0$:
\begin{equation*}
\norm f. \operatorname {exp}(L ^{\alpha }). \simeq 
\sup _{p>1} p ^{-1/\alpha } \norm f. p. \,. 
\end{equation*}
\end{proposition}

Comparing this Proposition to the results of the previous section, we see that 
an equivalent form of the Khintchine Inequalities for Rademacher random variables 
is as follows. 

\begin{theorem}\label{t.kexp} For 
all square summable coefficients $ c_j$ we have the following estimate 
for the Rademacher sequence of random variables $ \{r_j\mid j\}$. 
\begin{equation*}
\NOrm \sum _{j} c_j r_j . \operatorname {exp}(L ^2 ). 
\lesssim  \norm \{c_j\} . \ell ^2 . \,. 
\end{equation*}
\end{theorem}

Likewise, the following reformulation of Theorem~\ref{t.martKhintchine} is frequently 
referred to as the Chang-Wilson-Wolff inequality. 

\begin{cww} We have this inequality valid for all Hilbert space valued martingales. 
\begin{equation*}
\norm f . \operatorname {exp}(L ^2 ). \lesssim \norm \operatorname S (f). \infty . \,. 
\end{equation*}
\end{cww}

For $ \alpha >0$, let $ \varphi _{\alpha } (x)$ be a symmetric convex function which 
equals $ \lvert  x\rvert (\log 3+\lvert  x\rvert ) ^{\alpha } $ for $ \lvert  x\rvert $ 
sufficiently large, depending upon $ \alpha $.\footnote{For $ \alpha \ge 1$, we can take this 
as the definition for all $ \lvert  x\rvert\ge 0 $.}
The Orlicz space $ L ^{\varphi _{\alpha }}$ 
is denoted as $ L ^{\varphi  _{\alpha }} = L (\log L) ^{\alpha }$.  
These are the spaces that we used in the statement of our Main Theorems.

The point of interest 
here is following Proposition, which is in the standard references we have cited. 

\begin{proposition}\label{p.dual} For $ 0<\alpha <\infty $,  the two 
Orlicz spaces $ L (\log L) ^{\alpha }$ and $ \operatorname {exp}(L ^{1/\alpha })$
are Banach spaces which are dual to one-another.  
\end{proposition}

In particular, the dual to $ L (\log L)$ is $ \operatorname {exp}(L)$, and 
the dual to $ \operatorname {exp}( L ^{2} )$ is $ L \sqrt {\log L}$.

\subsection*{Discrepancy}

We recall some definitions and facts about Discrepancy 
which  are well represented in the literature, 
\cites{MR0066435,MR554923,MR903025}.

 Each dyadic interval has a left and right half, $ {I_{\textup{left}}}, {I_{\textup{right}}}$ 
respectively,  which are also dyadic.  Define the
 Haar function associated with $ I$ by  
 \begin{equation*}
h_I \coloneqq -\mathbf 1 _{I_{\textup{left}}}+ \mathbf 1 _{I_{\textup{right}}}
\end{equation*}
 Note that this is an $ L ^{\infty }$ normalization of these functions.

 In dimension $d $, a \emph{dyadic rectangle} is a product of dyadic intervals, thus an 
element of 
 $\mathcal D^d $.   A Haar function associated to 
 $R $ we take to the be product of the Haar functions associated 
 with each side of $R $, namely  for $ R_1\times \cdots \times  R_d$, 
 \begin{equation*}
 h_{R_1\times \cdots \times  R_d }(x_1 ,\dotsc, x_d) {} \coloneqq \prod _{t=1} ^{d}h _{R_t}(x_t)\,.
 \end{equation*}
We will concentrate on rectangles of a fixed volume, contained in 
  in $[0,1]^d $.

We call a function $f$ an \emph{$\mathsf r$ function  with parameter 
$ \vec r= (r_1 ,\dotsc, r_d)$} if $ \vec r\in \mathbb N ^{d}$, and 
\begin{equation}
\label{e.rfunction} f=\sum_{ R \in \mathcal R _{\vec r}}
\varepsilon_R\, h_R\,,\qquad \varepsilon_R\in \{\pm1\}\,. 
\end{equation}
In this sum, we let 
\begin{equation*}
\mathcal R _{\vec r} \coloneqq \{ R= (R_1 ,\dotsc, R_d) \mid R \ \textup{dyadic}\,,\ 
\lvert  R_t\rvert= 2 ^{- r_t}\,, \ 1\le t \le d \}\,. 
\end{equation*}
We will use $f _{\vec r} $ to denote a generic $\mathsf r$ function.   
A fact used without further comment is that $ f _{\vec r} ^2 \equiv 1$. 

Let $ \abs{ \vec r}= \sum _{t=1} ^{d} r_t=n$, which we refer to the 
index of the $ \mathsf r$ function. And let
\begin{equation*}
\mathbb H _n ^{d} \coloneqq \{\vec r \in \{1 ,\dotsc, n\} ^d 
 \mid \abs{ \vec r}=n\}\,. 
\end{equation*}
It is fundamental to the subject that $\mathbb H _n ^{d}  \simeq n ^{d-1}$. 
We refer to  $ \{f _{\vec r} \mid r\in \mathbb H _N ^{d}\}$ as hyperbolic $ \mathsf r$ functions.

The next four Propositions are standard.

\begin{proposition}\label{p.rvec} 
For any selection  $ \mathcal A_N$ of $ N$ points in the unit cube the following 
holds. 
Fix $ n$ with $ 2N< 2 ^{n}\le 4N$.  
For each $\vec r\in \mathbb H_{n} ^{d}$, there is an $\mathsf r$ function $f _{\vec r} $ with 
\begin{equation*}
\ip D_N, f _{\vec r}, \gtrsim 1\,. 
\end{equation*}
\end{proposition}

 \begin{proof}
 There is a very elementary one dimensional fact: For all dyadic intervals $I$, 
 \begin{equation} \label{e.veryelementary}
 \int _{0}  ^{1 } x   \cdot   h _{I}(x) \; dx =\tfrac 14 \abs{ I} ^2  \,.
 \end{equation}
 This immediately implies that in any dimension 
 \begin{equation} \label{e.vv}
 \ip   \abs{ [0,\vec x)} , h_R(\vec x), = 4 ^ {-d}\abs{ R} ^2  \,.
 \end{equation}

   Recall that $\mathcal A_N$, the distribution of 
 $N$ points in the unit cube, is fixed.  Call a cube $R\in \mathcal R _{\vec r}$ 
 \emph{good} if $R$ does {\bf not} intersect $\mathcal A_N$, otherwise call it \emph{bad}. 
  Set 
  \begin{equation}  \label{e.f_r}
  f _{\vec r} {} \coloneqq \sum _{R\in \mathcal R _{\vec r}}
  \operatorname {sgn} (\ip D_N, h_R,) h_R \,. 
  \end{equation}

  Each bad rectangle contains at least one point in $\mathcal A_N$, and  $2^ n\ge2N$, so 
  there are at least $N$ good rectangles.  Moreover, one should observe that the counting function 
 $\sharp (\mathcal A_N\cap [0,\vec x))$ is orthogonal to $ h_R$, for each good rectangle $
 R$.   Thus, by \eqref{e.vv},  
 \begin{equation*}
 \ip D_N, h_R,=-N\ip   \abs{ [0,\vec x)} , h_R(\vec x),=N2 ^ {-2n-2d}
\lesssim - 2  ^{-n}\,. 
 \end{equation*}
 Hence, we can estimate 
 \begin{align*}
 \ip D_N, f _{\vec r}, \ge 
 \sum _{\substack{R\in \mathcal R _{\vec r}\\ \text {$R$ is good} } }
 \ip D_N, h_R, \gtrsim 2^  {-n} \sharp \{
 R\in \mathcal R _{\vec r}\mid \text {$R$ is good} \} 
 \gtrsim 1\,. 
 \end{align*}
 And so our proof is complete.

 \end{proof}

\begin{proposition}\label{p.>n} Let $ f _{\vec s}$ be any $ \mathsf r$ function 
with $ \abs{ \vec s}>n$.  We have 
\begin{equation*}
\abs{ \ip D_N, f _{\vec s}, } \lesssim N {2 ^{-\abs{ \vec s}}}\,. 
\end{equation*}

\end{proposition}

 \begin{proof}
 This is a brute force proof.  Consider the linear part of the Discrepancy function. 
 By (\ref{e.veryelementary}), we have 
 \begin{equation*}
 \abs{ \ip N \prod _{j=1} ^{d} x_j, f _{\vec s}, }\lesssim N 2 ^{-\abs{ \vec s}} \,,
 \end{equation*}
 as claimed. 
 
 Consider the part of the Discrepancy function that arises from the point set.  
 Observe that for any point $\vec x_0$ in the point set, we have 
 \begin{equation*}
 \abs{ \ip \mathbf 1 _{[\vec 0, \vec x_0)} , f _{\vec s}, } \lesssim 2 ^{- \abs{ \vec
 s}}\,.
 \end{equation*}
 Indeed, of the different Haar functions that contribute to $ f _{\vec s}$, there 
 is at most one with non zero inner product with the function 
 $\mathbf 1 _{[\vec 0, \vec x)} (\vec x_0) $ as a function of $ \vec x$.  It is the 
 one rectangle which contains $ x_0$ in its interior. Thus the inequality above follows. 
 Summing it over the $ N$ points in the point set finish the proof of the Proposition. 
 
 \end{proof}

In two dimensions,  the decisive product rule holds.  It is as follows, and we omit the
proof.  

\begin{proposition}\label{p.productRule}  
In dimension $ d=2$ the following holds. 
Let $ \vec r_1 ,\dotsc, \vec r_k$ be   elements of $ \mathbb H _n ^{2}$ 
where one of the vectors occurs an odd number of times. 
Then, the product $ \prod _{j=1} ^{k} f _{\vec r}$ is also a $ \mathsf r$ function. 
If the $ \vec r_j$ are distinct and $ k\ge 2$, the product has index larger than $ n$. 
\end{proposition}

\begin{proposition}\label{p.Not2Many} 
In dimension $ d=2$ the following holds. 
Fix a collection of $ \mathsf r$ functions 
$\{ f _{\vec r} \mid \vec r \in \mathbb H_{n} ^{2} \}$. 
Fix an integer $ 2\le v \le n$ and $ \vec s $ with $ \abs{ \vec s}\ge n+ v -1$. 
Let $ \operatorname {Count} (\vec s ; v) $ be 
{the number of ways to  choose distinct 
$ r_1 ,\dotsc, r_v\in \mathbb H _n ^2 $  so that $ \prod _{w=1} ^{v} f _{\vec r_w}$
is an $ \vec s$ function.}
We have 
\begin{equation} \label{e.Not2Many}
\operatorname {Count} (\vec s ; v) 
=  {\abs {\vec s}-n-1 \choose v -2 }\,. 
\end{equation}
\end{proposition}

 \begin{proof}
 Fix a vector $ \vec s$ with $ \abs{ \vec s}>n$, and suppose that 
 \begin{equation*}
 \prod _{w=1} ^{v} f _{\vec r_w} 
 \end{equation*}
 is an $ \vec s$ function.  Then, the maximum of the first coordinates of the $ \vec r_w$
 must be $ s_1$, and similarly for the second coordinate.  Thus, the vector $ s$ 
 completely specifies two of the $ \vec r_w$.  
 
 The remaining $ v-2$ vectors must be distinct, and take values in the first 
 coordinate that are less than $ s_1$, and in the second less than $ s_2$. 
 The hyperbolic assumption then says that there are at most $ \abs{ \vec s}-n-1$ 
 possible choices for these vectors, and from them we select $ v-2$.  This completes 
 the proof.

 \end{proof}

The next Lemma is elementary, and has probably been observed before, but 
is not a standard fact. 
For integers $ 1\le v \le n$, let 
\begin{equation}\label{e.G}
G_v \coloneqq
\sum _{\{\vec r_1 ,\dotsc, \vec r_v \}\subset \mathbb H _n ^2 \\
} \prod _{w=1} ^{v} f _{\vec r_w}\,.
\end{equation}
This sum is over all  subsets of $ \mathbb H _n ^2 $ of cardinality $ v$, that is 
the vectors $ \{\vec r_1 ,\dotsc, \vec r_v \}$ are all distinct. 

\begin{lemma}\label{l.G} For odd integers $ k$ we have 
\begin{gather} \label{e.G=}
\Biggl[ \sum _{\vec r\in \mathbb H _n ^2 } f _{\vec r} \Biggr] ^{k} 
= 
\sum _{\substack{v=1\\ \textup{$ v$ odd} }} ^{\min (k,n)} 
\gamma (k,n) \cdot  n ^{(k-v)/2} \cdot 
G_v \,, 
\\  \label{e.gamma}
0<\gamma (k,n) \le 
\frac {k !} { (k-v)!} [ C_0  (k-v) ] ^{(k-v)/2} \le k! C ^{ k-v}\,. 
\end{gather}
where $ C _0 $ and $ C$ are  absolute constants. 
\end{lemma}

\begin{proof}
The  point is that $ f _{\vec r} ^2 \equiv 1$, and that $ k$ is odd. 
Thus, even products 
of distinct $ \mathsf r$ functions cannot arise in this product. 
(There would be a different formulation for $ k $ even.)

For a fixed subset $ \{ \vec r_1 ,\dotsc, \vec r_v\} \subset \mathbb H _n ^2 $, 
with $ v$ odd and at most $ k$,  we consider the number of ways that the product 
\begin{equation*}
\prod _{w=1} ^{v} f _{\vec r_w}
\end{equation*}
 can arise in the left hand side of (\ref{e.G=}).
From the  $ k$ terms  on the left in (\ref{e.G=}), 
we choose $ v $  from which we take one of the pre-specified $\mathsf r $ 
functions $ f _{\vec r_1} ,\dotsc, f _{\vec r_v}$ in some order. 
There are $ k!/ (k-v!)$ ways to do this. 

In the remaining $ k-v$ terms, an even number of terms, we need to estimate $ \gamma ' (k-v)$, 
the number  of ways to assign $ \mathsf r$ functions to these, so that they can be divided into distinct
pairs of equal $ \mathsf r$ functions.  This is a somewhat tricky combinatorial 
problem.  Let us observe that this problem arises in a concrete way in 
one approach to the Khintchine inequalities.

Let $ \varphi _{\vec r}$ be the $ \mathsf r$ function as in \eqref{e.rfunction}, when 
all the choices of signs $ \varepsilon _{R}$ are identically one.  Thus, $ \varphi _R$ are 
Rademacher functions.  Then, the term we are seeking to bound is 
\begin{equation*}
\gamma ' (k-v)=\int _{[0,1] ^{2}} \Biggl[\sum _{\lvert  \vec r\rvert=n } \varphi _{\vec r} \Biggr] ^{k-v} \; dx
\end{equation*}
Indeed, expanding the $ (k-v)$th power on the right, the integral of the product below 
\begin{equation*}
\int _{[0,1] ^{2}} \prod _{u=1} ^{k-v} \varphi _{\vec r _u}\; dx
\end{equation*}
will be non-zero iff the $ \{\vec r_1 ,\dotsc, \vec r _{k-v}\}$ can be written as a 
disjoint union of pairs of equal vectors.  

Classical proofs of the Khintchine inequality, see \cite{MR0500056}, estimate this norm directly. 
We can use the well-known best constants to estimate as below.  The best constants 
in the range we are interested were established in \cite{MR0438089}.  The only information 
that we use here is the asymptotic order, which follows from other sources, such as 
\cite{MR800004,MR1018577,MR850744}, as well as  references therein.  For an absolute constant $ C$ we
have 
\begin{align*}
\int_{[0,1] ^{2}}  \Biggl[\sum _{\lvert  \vec r\rvert=n } \varphi _{\vec r} \Biggr] ^{k-v} \; dx
 = 
\NOrm \sum _{\lvert  \vec r\rvert=n } \varphi _{\vec r} . k-v . ^{k-v} 
 \le [ C \sqrt {k-v} \cdot \sqrt n ] ^{k-v}\,. 
\end{align*}
This completes our proof. 

\end{proof}

\section{The Proof of Hal{\'a}sz' Theorem} 

The point of this section is to provide a new proof of Hal{\'a}sz' Theorem, 
on the $ L ^{1}$ norm of the Discrepancy function in two dimensions. 
It suffices to prove this for sufficiently large $ N$. 
The proof is by way of duality.
Fix the point distribution $ \mathcal A_N
\subset [0,1] ^2 $.   Set $ 2N<2 ^{n} \le 4N$, so that $ n \simeq \log N$. 
Proposition~\ref{p.rvec} provides us with  
$ \mathsf  r$ functions $ f _{\vec r} $ for $ r\in \mathbb H _n ^{2}$. 
We show that $ \ip D_N, \Psi , \gtrsim \sqrt n$, where 
\begin{equation} \label{e.Psi}
\Psi \coloneqq \sin \Bigl( \epsilon \cdot   n ^{-1/2} \sum _{\vec r\in \mathbb H _n ^2 }  f _{\vec r}\Bigr)\,.
\end{equation}
Here $ \epsilon >0$ is a small positive constant.  It needs to be smaller than the 
$ 1/C_1$ for $ C_1$ as in \eqref{e.delta}. 
Of course $ \Psi $ is a bounded function, so  this will prove the Theorem.  

For any fixed $ n$, $  n ^{-1/2} \sum _{\vec r\in \mathbb H _n ^2 }  f _{\vec r}$ is 
a bounded function, hence the Taylor series expansion of $ \Psi $ is absolutely convergent. 
Lemma~\ref{l.G} shows that the  Taylor series for $ \Psi $ satisfies 
\begin{align} \nonumber 
\Psi & = 
\sum _{ \textup{ $ k$ odd} } 
\frac { (-1) ^{(k+1)/2}} {k!}  n ^{-k/2} \epsilon  ^{k}
\Bigl[ \sum _{\vec r\in \mathbb H _n} f _{\vec r} \Bigr] ^{k} 
\\ \nonumber 
& = 
\sum _{ \textup{ $ k$ odd} } 
 { (-1) ^{(k+1)/2}}   \epsilon ^{k} \sum _{\substack{ v=1\\ \textup{ $ v$ odd} }} ^{n} 
\frac {\gamma (k,v)} {k!} n ^{-v/2} G_v 
\\
& =  \label{e.finalexpand}  
\sum _{\substack{ v=1\\ \textup{ $ v$ odd} }} ^{n}   \delta (v) 
n ^{-v/2} G_v \,.  
\end{align}
In this last display, we can take $ \delta (v)$ to be a sequence of constants with 
$ \delta (1) = C_0\epsilon $, for a choice of constant  $ C _0$, which depends only on    
 the constant $ C$ in \eqref{e.gamma}.   For $ v >1$, 
\begin{equation}\label{e.delta}
\lvert  \delta (v)\rvert \le (C_1\epsilon)  ^{v}  \,, \qquad v>1\,. 
\end{equation}
for an   absolute constant $ C_1$.

We turn our attention to the terms in (\ref{e.finalexpand}). 
Now, by construction, we have 
\begin{equation} \label{e.ll}
\ip D_N ,  \delta _0n ^{-1/2} G_1 , \gtrsim  \delta _0 n ^{-1/2} \sum _{\vec r\in \mathbb H _n ^2 } 
\ip D_N, f _{\vec r}, \gtrsim  \delta _0 n ^{1/2} \simeq \epsilon  \sqrt {\log N}\,. 
\end{equation}

As for the terms $ 3 \le v \le n$, note that by Proposition~\ref{p.>n}, 
Proposition~\ref{p.Not2Many} and the definition of $ G_v$, we have 
\begin{align*}
\abs{ \ip D_N, G_v , } 
& \lesssim  N \sum _{s=n+v-1} ^{2n}  2^{-s} {s-n-1 \choose v-2}\,. 
\end{align*}
And so we  estimate as follows, using \eqref{e.finalexpand} and \eqref{e.delta}. 
	Here is convenient that the 
sum is only over odd $v \ge 3 $. 
\begin{align*}
\sum _{\substack{ v=3\\ \textup{$ v$ odd } }} ^{n} 
 (C_1 \epsilon ) ^{v} n ^{-v/2} \abs{ \ip D_N, G_v , } 
 & 
 \lesssim 
N  \sum _{\substack{ v=3\\ \textup{$ v$ odd } }} ^{n} 
\sum _{s=n+v-1} ^{2n} (C_1 \epsilon ) ^{v} 2 ^{-s} n ^{-v/2} {s-n-1 \choose v-2}
\\
& \lesssim N  n ^{-1} 
\sum _{s = n+3} ^{2n} 2 ^{-s} \sum _{v=0} ^{ s-n-1 }  (C_1 \epsilon ) ^{v}
n ^{-v/2} {s-n-1 \choose v} 
\\
& \lesssim  n ^{-1} \sum _{s=0} ^{n }  2 ^{-s}
(1+ (C_1 \epsilon ) ^2 n ^{-1/2} )^{s} 
\\
& \lesssim   n ^{-1}\,. 
\end{align*}
This estimate holds for $ n$ sufficiently large, and 
tends to zero with $ n$,  so 
we can combine it with \eqref{e.ll} to 
 prove the  Hal{\'a}sz Theorem for sufficiently 
large $ N$.

\section{The Proof of Theorem~\ref{t.beyondHalasz}} 
 
The proof is by duality.  The Orlicz space in which we seek an estimate of $ D_N  $
is $ L ^{1} (\log L) ^{ (d-2)/2}$.  The dual to this space, by Proposition~\ref{p.dual}, 
is the  exponential Orlicz 
space 
$ \operatorname {exp} (L ^{2/ (d-1)})$.  Namely we define an $ \Phi $ so that 
\begin{equation*}
\ip D_N, \Phi , \gtrsim (\log N) ^{ (d-1)/2}\,, 
\qquad 
\norm \Phi . \operatorname {exp} (L ^{2/ (d-1)}). \lesssim 1\,. 
\end{equation*}

Given the point 
set $ \mathcal A_N\subset [0,1] ^{d}$, we fix $ 2N\le 2 ^{n}<4N$, 
and take the $ \mathsf r$ functions as in Proposition~\ref{p.rvec}. 
We construct our test function as follows.

 For $ \vec s\in \mathbb N ^{d-2}$, let 
\begin{equation*}
F _{\vec s} \coloneqq \sum _{\substack{\vec r \in \mathbb H_n ^{d} 
	\\ r_j = s_j\,, \ 1\le j\le d-2}} f _{\vec r}
\end{equation*}
That is, we sum over vectors  $ \vec r$ which equal $ \vec s$ in the first 
$ d-2$ coordinates.   Clearly, we need to require at a minimum that 
$ \lvert  \vec s\rvert< n $, so that there are only $ \le  n ^{d-2}$ possible values 
$ \vec s$ for which the definition above is non-zero. 
In dimension $ d=2$  we interpret this as no restriction on the coordinates. 
In dimension $ d=3$ only the  first coordinate is restricted. 
 
 Our test function is 
 \begin{equation*}
\Phi \coloneqq n ^{- (d-2)/2} \sum _{\lvert  \vec s\rvert\le \tfrac 34 n } \sin ( \epsilon  n
^{-1/2} F _{\vec s} )\,. 
\end{equation*}
Here $ \epsilon >0$ is the small positive constant of the previous section. 
The argument in the two dimensional case can be modified to show that 
\begin{equation*}
\ip D_N, \sin ( n ^{-1/2} F _{\vec s} ) , \gtrsim n ^{1/2} \,, 
\qquad \lvert  \vec s\rvert\le \tfrac34 n\,.  
\end{equation*}
Therefore, it follows that $ \ip D_N, \Phi , \gtrsim n ^{(d-1)/2}$.

Let us check the integrability properties of $ \Phi $.  
That is, using Proposition~\ref{p.comparable}, we should verify that 
\begin{equation} \label{e.check}
\norm \Phi .p. \lesssim p ^{ (d-2)/2}  \,. 
\end{equation}
We shall find that $ d-2$ applications of the Littlewood--Paley 
are enough to prove this estimate. 

Consider the function $ \Phi $ with the variables $ x_2 ,\dotsc, x_d$ fixed.  
The terms below  form a martingale difference sequence in $ \sigma _1$:  
\begin{equation*}
\sum _{\substack{\vec s\\ s_1= \sigma _1}} 
 n ^{- (d-2)/2}  \sin ( \epsilon  n ^{-1/2} F _{\vec s} )\,. 
\end{equation*}
The square function of this martingale, which we denote by $ \operatorname S_ 1 (\Phi )$ is 
\begin{equation*}
\operatorname S _1 (\Phi ) 
=n ^{- (d-2)/2}  \Biggl[ \sum _{\sigma _1}  \ABs{ 
\sum _{\substack{\vec s\\ s_1= \sigma _1}} 
 \sin ( \epsilon  n ^{-1/2} F _{\vec s} )} ^2 \Biggr] ^{1/2} \,. 
\end{equation*}
A useful (but not absolutely essential) point to observe is that the 
term on the right above is a Hilbert space valued function 
\begin{equation*}
\operatorname S _1 (\Phi ) 
= n ^{- (d-2)/2} 
\NORm \Biggl\{\sum _{\substack{\vec s \in \mathbb N ^{d-2}\\ s_1= \sigma _1}} 
 \sin ( \epsilon  n ^{-1/2} F _{\vec s} )  \mid \sigma _1\in \mathbb N  \Biggr\}. \ell ^{2} (\sigma _1) .
\end{equation*}
Continuing by induction, holding the variables $ x_3 ,\dotsc, x_d$ constant above, 
we see that the expression above is a Hilbert space martingale, with martingale difference 
sequence 
\begin{equation*}
\Biggl\{  \sum _{\substack{\vec s \in \mathbb N ^{d-2}\\ s_1= \sigma _1 \,, \, s_2=\sigma _2 }} 
 \sin ( \epsilon  n ^{-1/2} F _{\vec s} )  \mid \sigma _1 \,,\, \sigma _2\in \mathbb N 
 \Biggr\}\,. 
\end{equation*}
We accordingly define 
\begin{equation*}
S_2 (\Phi ) 
= \Biggl[ \sum _{\sigma _1} \sum _{\sigma _2}  \ABs{ 
\sum _{\substack{\vec s\\ s_1= \sigma _1 \,,\,  s_2= \sigma _2}} 
 \sin ( \epsilon  n ^{-1/2} F _{\vec s} )} ^2 \Biggr] ^{1/2} \,.  
\end{equation*}
This notation is extended to $ d-2$, where 
\begin{equation}  \label{e.icww}
S_ {d-2} (\Phi ) 
= n ^{- (d-2)/2} 
\Biggl[   
\sum _{\substack{\vec s}} 
 \sin  ^2 ( \epsilon  n ^{-1/2} F _{\vec s} )  \Biggr] ^{1/2}
  \le  1\,. 
\end{equation}

We now apply the Littlewood-Paley inequalities, as phrased in Theorem~\ref{t.martKhintchine}, 
to conclude that  for $ 2\le p <\infty $, 
\begin{align*}
\norm \Phi .p. & \le C \sqrt p \norm S_1 (\Phi ).p.
\\
& \le C ^2 p \norm S_2 (\Phi ).p. 
\\
&\ \ \vdots
\\
&\le C ^{p-2} p ^{(d-2)/2} \norm S _{d-2} (\Phi ). p. \lesssim p ^{(d-2)/2} \,. 
\end{align*}
That is, \eqref{e.check} holds, finishing our proof.   Our proof of the main result is 
complete. 

We comment, that with the $ L ^{\infty } $ bound on the iterated square function, especially
in \eqref{e.icww}, we are employing a variant of the Chang-Wilson-Wolff inequality 
as described in \cite{MR850744,MR1439553}.

\smallskip 

In dimension $ d=3$, an outstanding conjecture is that there is the universal 
estimate 
\begin{equation*}
\norm D_N.1. \gtrsim \log N\,,  
\end{equation*}
valid for all point sets $ \mathcal A_N$.  
To attempt to prove this estimate, 
it would be natural to consider the function 
\begin{equation*}
\Psi \coloneqq 
\sin \Bigl(  n ^{-1/2}\sum _{a} \sin \Bigl(  n ^{-1/2}\sum _{\substack{ 
\vec r\in \mathbb H _N ^{3}\\ r_1=a }} f _{\vec r} \Bigr) \Bigr)\,. 
\end{equation*}
But there are difficulties to this line of approach that are 
similar to those that complicate 
the argument in J.~Beck's paper \cite{MR1032337}. 
Indeed, \emph{any} improvement in our main Theorem, in dimension $ d=3$ 
appears to be of interest. 

\section{A Remark on Hardy Spaces and the Discrepancy Function} 

For functions $ f$ on $ \mathbb R ^{d}$, and  $ 1\le t\le d$ we 
define the integration operator in coordinate $ t$ by 
\begin{equation}\label{e.Int}
\operatorname {Int} _{t} f (\vec x) \coloneqq \int _{0} ^{1} 
f (x_1 ,\dotsc, x_d) \; d x _{t}\,. 
\end{equation}
For a point set $ \mathcal A_N$ of cardinality $ N$, with Discrepancy function 
$ D_N$, we define 
\begin{equation}\label{e.tilde}
\widetilde D_N = \Bigl(\operatorname {Id}- \sum _{T\subset \{1 ,\dotsc, d\}} 
(-1) ^{\lvert  T\rvert }
\prod _{t\in T} \operatorname {Int}_t\Bigr) D_N\,. 
\end{equation}
Here $ \operatorname {Id}$ is the Identity operator. We prove 

\begin{theorem}\label{t.hardyD} For $ d\ge 2$ and all $ 0<p\le 1$, we have the estimate 
\begin{equation*}
\Norm \widetilde D_N . H ^{p}. \gtrsim (\log N) ^{(d-1)/2}\,. 
\end{equation*}
Here, $ H^p$ is the dyadic  real-valued Hardy space of $ d$ parameters on $ \mathbb R ^{d}$. 
\end{theorem}

The Hardy space we have in mind in this Theorem is given in the next Definition. 

\begin{definition}\label{d.hardy} The 
dyadic  real-valued Hardy space of $ d$ parameters on $ \mathbb R ^{d}$ 
is the set of functions $ f \;:\;  \mathbb R ^{d} \longrightarrow \mathbb R $ such that 
the following norm is finite. 
\begin{align} \label{e.Hpdef}
\norm f . H ^{p}. & \coloneqq  \norm \operatorname M f .p. \,,  
\\
\label{e.Mdef}
\operatorname M f (x) &\coloneqq \sup _{\substack{x\in R\\ R \ \textup{dyadic} }}
\frac 1 {\lvert  R\rvert } \int _{R} f (y)\; dy \,. 
\end{align}
Here, at a given point $ x$, the supremum is over all dyadic rectangles in $ \mathbb R ^{d}$
that contain $ x$. Note that we do not take the absolute value of the function $ f$. 
\end{definition}

The Hardy space in question, is a `real-variable' extension of the classical space 
of analytic functions introduced by G.~H.~Hardy.  It is a subtle object, 
whose crucial properties have been 
identified by Alice Chang, Robert Fefferman, 
\cite{MR584078,cf1}, 
Jean-Lin Journ\'e \cite{MR826486} and Jill Pipher \cite{MR850744}, among 
others.  The dyadic setting, of interest to us, is specifically addressed in the 
paper \cite{MR539351}. 
See especially the paper \cites{cf1}, as well as the references in these papers. 

There is an alternative definition, which we recall here. 
Let us define the \emph{Square Function} by 
\begin{equation*}
\operatorname S (f) \coloneqq \Biggl[ \sum _{R} \frac {\lvert  \ip f, h_R, \rvert ^2  } 
{\lvert  R\rvert ^2  } \mathbf 1_{R}\Biggr] ^{1/2} \,. 
\end{equation*}
The sum is over dyadic rectangles in $ \mathbb R ^{d}$.

\begin{Hardy}  The following equivalence of norms holds.  
\begin{equation*}
\norm \operatorname M f . p. \simeq \norm \operatorname S (f). p. \,, 
\qquad 0<p \le  1\,. 
\end{equation*}
\end{Hardy}

It is well-known that Hardy spaces are an appropriate 
substitute for $ L ^{p}$ spaces for a variety of issues concerning Harmonic Analysis. 
In particular, there is a third equivalent definition of the Hardy space norm in terms of 
Hilbert transforms, which we omit in this discussion.  There is a fourth `atomic
decomposition' approach that we also omit (See \cite{atomic}), referring all these issues to the cited
references.

\bigskip

The point of the subtraction in \eqref{e.tilde} is that $ D_N$ is not \emph{a priori} 
a member of Hardy space, but $ \widetilde D_N$ is.  Note that the conclusion of the 
Theorem provides partial support for the Conjecture below, which is an extension of the 
$ L ^{1}$ Conjecture mentioned in the introduction.  There is nothing known about this 
Conjecture, even in dimension $ d=2$. 

\begin{Lp}
In dimension $ d\ge 2$ we have the estimate 
\begin{equation*}
\norm D_N .p. \gtrsim  (\log N ) ^{ (d-1)/2}\,, 
\qquad 0<p<1\,. 
\end{equation*}
\end{Lp}

\begin{proof}[Proof of Theorem~\ref{t.hardyD}.] 
We use the following elementary fact.  Let $ G_1 ,\dotsc, G_J$ be 
measurable subsets of a probability space $ (\Omega ,\mathbb P )$, with $ \mathbb P (G_j)\ge 
\tfrac 12$ for all $ 1\le j \le J$.  Then, 
\begin{equation*}
\mathbb P \Bigl(\sum _{j= 1} ^{J} \mathbf 1 _{G_j} > J/4 \Bigr)\ge \tfrac 14\,. 
\end{equation*}
Indeed, suppose this is not the case, then we have the contradiction 
\begin{align*}
\frac J 2&\le \NOrm \sum _{j= 1} ^{J} \mathbf 1 _{G_j} .1. 
\\
& \le J/4+ J \cdot \mathbb P \Bigl(\sum _{j= 1} ^{J} \mathbf 1 _{G_j} > J/4 \Bigr)
< \frac J 2\,. 
\end{align*}
It follows that we have 
\begin{equation}\label{e.p-norm}
\NOrm \sum _{j= 1} ^{J} \mathbf 1 _{G_j}  .p. \gtrsim J\,, \qquad 0<p \le \infty \,. 
\end{equation}

Fix a point set $ \mathcal A_N$, and let $ 2N\le 2 ^{n}<4N$. 
For a vector $ \lvert  \vec r\rvert=n $, let 
\begin{equation*}
\mathcal G _{\vec r} \coloneqq  \{R \;:\;   \lvert  R_t\rvert= 2 ^{r_t}\,,\ 1\le t \le d \,,\ 
R \cap \mathcal A_N=\emptyset\}\,. 
\end{equation*}
These rectangles, which avoid the point set $ \mathcal A_N$  are the `good rectangles' in  
the proof of Proposition~\ref{p.rvec}.  Set $ G _{\vec r}=\bigcup _{R\in \mathcal G _{\vec
r}} R$. 
By the pigeonhole principle, $ \lvert  G _{\vec r}\rvert\ge \tfrac 12 $.  Moreover, 
a standard computation, see \eqref{e.vv},  shows that 
\begin{equation*}
\frac {\lvert  \ip \widetilde D_N, h _{R},\rvert ^2  } {\lvert  R\rvert ^2  } \gtrsim 1\,, 
\qquad R\in \mathcal G _{\vec r}\,. 
\end{equation*}
The implied constant only depends upon dimension.

We conclude a lower bound on the Square Function of $ \widetilde D_N$:
\begin{equation*}
\operatorname S (\widetilde D_N) \gtrsim  \Biggl[ \sum _{ \lvert  \vec r\rvert=n } 
\mathbf 1 _{G _{\vec r}}\Biggr] ^{1/2} \,. 
\end{equation*}
As there are $ n ^{d-1} \simeq (\log N) ^{d-1}$ choices of $ \vec r$, we conclude 
from \eqref{e.p-norm} that we have 
\begin{equation*}
\norm \operatorname S (\widetilde D_N) .p. \gtrsim n ^{(d-1)/2}\,, 
\qquad 0<p\le 1\,. 
\end{equation*}
By the definition of the Hardy space norm, this completes our proof. 
\end{proof}

 \end{document}